\documentclass[draft]{amsart}
\usepackage{amssymb}
\theoremstyle{plain}
\newtheorem{teo}{Theorem}
\newtheorem{cor}{Corollary}

\newcommand{\cl}{{\rm cl\,}}

\begin{document}
\date{}
\title{A note on lower bounds of martingale measure densities}
\author{Dmitry Rokhlin}
\address{Dmitry Rokhlin, Department of Mechanics and
         Mathematics, Rostov State University, Zorge str.~5,
         Rostov-on-Don, 344090, Russian Federation}
\email{rokhlin@math.rsu.ru}
\author{Walter Schachermayer}
\address{Walter Schachermayer, Vienna University of Technology,
Wiedner Hauptstrasse 8-10/105, A-1040 Wien, Austria and
Universit\'{e} Paris Dauphine, Place du Mar\'{e}chal de Lattre de Tassigny,
F-75775 Paris Cedex 16, France}
\email{wschach@fam.tuwien.ac.at}
\thanks{The second author gratefully acknowledge financial support
from the Austrian Science Fund (FWF) under Grant P15889 and from
Vienna Science and Technology Fund (WWTF) under Grant MA13}
\dedicatory{Dedicated to the memory of Joe Doob}
\begin{abstract}
For a given element $f\in L^1$ and a convex cone $C\subset L^\infty$,
$C\cap L^\infty_+=\{0\}$ we give necessary and sufficient
conditions for the existence of an element $g\ge f$ lying in
the polar of $C$. This polar is taken in $(L^\infty)^*$ and
in $L^1$. In the context of mathematical finance the main result
concerns the existence of martingale measures, whose densities are
bounded from below by prescribed random variable.
\end{abstract}
\maketitle
\section{Introduction}
\setcounter{equation}{0}
Let $(\Omega,\mathcal F,\mathbf P)$ be a probability space.
Consider a convex cone $C\subset L^\infty=
L^\infty(\Omega,\mathcal F,\mathbf P)$, satisfying the condition
\begin{equation}
C\cap L^\infty_+=\{0\},
\end{equation}
where $L^\infty_+$ is the non-negative orthant of $L^\infty$.
Typically, $C$ consists of random variables,
dominated by stochastic integrals $\int_0^T H_t\,dS_t$
(compare \cite{DS05}). Here
$S=(S_t)_{0\le t\le T}$ is a semimartingale, describing
the stock-price process and $H=(H_t)_{0\le t\le T}$ is a predictable
$S$-integrable process, belonging to some class of admissible
trading strategies. Assumption (1.1) is usually referred to as
the no-arbitrage condition. Note, that the cases of
transaction costs, portfolio constraints and infinitely many
assets can also be incorporated in this framework.

Furthemore, consider the polar of $C$, taken in $L^1=
L^1(\Omega,\mathcal F,\mathbf P)$:
\begin{equation}
 \{y\in L^1:\int_\Omega xy\,d\mathbf P\le 0,\ x\in C\}.
\end{equation}
For the case of a bounded process $S$,
the set (1.2) is generated by densities of absolutely
continuous martingale measures.
In this note we discuss the following question:
\begin{itemize}
\item[\bf(Q)] Let $f\in L^1$. Under what conditions there
exists an element $g\in L^1$ in the polar of $C$ such that
$g\ge f$?
\end{itemize}

In fact, this question concerns the existence of a martingale
measure $\mathbf Q$, whose density is bounded from below by
the prescribed random variable $f$ up to a multiplicative
constant $\alpha>0$: $d\mathbf Q/d\mathbf P\ge \alpha f$.

Sometimes it is usefull to take the polar of $C$
in $(L^\infty)^*$, the dual space of $L^\infty$: see, e.g.
\cite{CSW01}. In our case it also
appears that an easier answer to the question
{\bf (Q)} can be given if $g$ is allowed to lie in
$(L^\infty)^*$: see Corollary 1 below and \cite{R04}.
The answer to this question in precise terms is given
in Corollary 2.

Our results are essentially the following.
Regard $f\in L^1$ as a functional on $L^\infty$, defined by the formula
$$ \langle x,f\rangle=\int_\Omega xf\,d\mathbf P. $$
Then the existence of the desired element $g$ is
equivalent to the boundness of $f$ from above on a certain subset
of the cone $C$. If $g$ is allowed to be an element of $(L^\infty)^*$,
this subset may be chosen as
$$ C_1=\{x\in C: x^-\le 1\ a.s.\}, $$
where $x^-=\max\{-x,0\}$.
If we seek for $g\in L^1$, such a subset should be somewhat bigger:
$$ C_V=\{x\in C: x^-\in V\}, $$
where $V$ is a neighbourhood of zero in the Mackey
topology $\tau (L^\infty,L^1)$.

\section{Answer to the question (Q)}
\setcounter{equation}{0}
We find it natural to examine the problem in a somewhat more
general context. Let $(X,\tau)$ be a locally convex-solid
Riesz space. It means that $X$ is a vector lattice, endowed
with a topology $\tau$, whose local base consists of
convex solid sets: see \cite{AB78} for details. For an element
$x\in X$, its positive part, negative part and absolute
value are denoted by $x^+$, $x^-$ and $|x|$. The set
$V\subset X$ is called solid if
the conditions $x\in V$, $|y|\le |x|$ imply that $y\in V$.

Consider a convex cone $C\subset X$, such that
\begin{equation}
C\cap X_+=\{0\},
\end{equation}
where $X_+=\{x\in X: x\ge 0\}$. Let
$V$ be a solid subset of $X$. Put
$$C_V=\{x\in C: x^-\in V\}.$$
Using the implication
\begin{equation}
x\le y\Longrightarrow x^-\ge y^-,
\end{equation}
it is elementary to check that
\begin{equation}
 C_V=C\cap (V+X_+).
\end{equation}

Denote by $X^*$ be the topological dual of $X$ with the order,
induced by the dual cone
$X^*_+=\{\xi\in X^*:\langle x,\xi\rangle\ge 0,\ x\in X_+\}$.
The polar of $C$ is taken in $X^*$:
$$ C^\circ=\{\xi\in X^*: \langle x,\xi\rangle\le 0,\ \ x\in C\}.$$
We use the customary notation  $\sigma(X^*,X)$ for the weak-star
topology and $|\sigma|(X,X^*)$ for the coarsest locally
convex-solid topology on $X$, compatible with the duality
$\langle X,X^*\rangle$ \cite{AB78}.
\begin{teo}
Let $(X,\tau)$ be a locally convex-solid Riesz space. Assume that
there exists a $\sigma(X^*,X)$-compact set $\Gamma\subset X^*_+$
such that the convex cone, generated by $\Gamma$ is
$\sigma(X^*,X)$-dense in $X^*_+$. Let $C\subset X$ be a convex
cone, satisfying (2.1).
Then for any $f\in X^*$ the following conditions are
equivalent:
\item[(i)] there exists a convex solid $\tau$-neighbourhood of zero
$V$ such that
$$\sup_{x\in C_V}\langle x,f\rangle<+\infty,\ \ \
   C_V=\{x\in C: x^-\in V\};$$
\item[(ii)] there exists $g\in C^\circ$ such that $g\ge f$.
\end{teo}

\begin{proof} (ii) $\Longrightarrow$ (i). Consider
the convex solid $|\sigma|(X,X^*)$-neighbourhood of zero
$$V=\{x\in X:\langle |x|,g-f\rangle\le 1\}.$$
Let $x\in C_V$, then
$$ \langle x,f\rangle=\langle x,g\rangle+\langle x,f-g\rangle
   \le \langle -x,g-f\rangle\le \langle x^-,g-f\rangle\le 1. $$

(i) $\Longrightarrow$ (ii). Let $\Gamma'$ be the
$\sigma (X^*,X)$-closed convex hull of the set $\Gamma\cup\{0\}$.
Consider the $\sigma (X^*,X)$-compact convex set
$$\Pi=(V-X_+)^\circ+\Gamma'=(V^\circ\cap X^*_+)+\Gamma'$$
and put
\begin{equation}
\lambda=\sup_{x\in C_V}\langle x,f\rangle.
\end{equation}

If the condition (ii) is false, we may apply the Hahn-Banach
theorem \cite[Chap.\,II, Th.\,9.2]{S66} to separate
the sets $f+\lambda\Pi$ and $C^\circ$ by an element $x\in X$:
$$ \sup_{\eta\in C^\circ}\langle x,\eta\rangle
   <\inf_{\zeta\in f+\lambda\Pi}\langle x,\zeta\rangle.$$
Since $C^\circ$ is a cone, we get $\langle x,\eta\rangle\le 0$,
$\eta\in C^\circ$. Thus, $x\in C^{\circ\circ}=\cl C$
by the bipolar theorem \cite[Chap.\,IV, Th.\,1.5]{S66},
where $\cl C$ is the closure of $C$ in any topology, compatible
with the duality $\langle X,X^*\rangle$, and
\begin{equation}
  \langle x,f\rangle+\lambda\inf_{\zeta\in\Pi}
  \langle x,\zeta\rangle>0.
\end{equation}

Furthemore, since $\inf_{\zeta\in\Pi}\langle x,\zeta\rangle\le 0$,
we conclude
that $\langle x,f\rangle>0$ and $x\not\in X_+$.
Indeed, for any $\tau$-neighbourhood of zero $W$ take an element
$y_W\in (\mu x+W\cap V)\cap C$, $\mu>0$. If $x^-=0$ then
$y_W\ge z_W$ for some $z_W\in V$. By (2.2) and the solidness of
$V$ we have $y_W^-\in V$. Thus, $\mu x\in\cl C_V$
for any $\mu>0$ and we obtain a contradiction,
since $\langle x,f\rangle>0$ and
$f$ must be bounded (from above) on $\cl C_V$.

Moreover, $\inf_{\zeta\in\Pi}\langle x,\zeta\rangle<0$, because
otherwise $x$ is non-negative on $\Gamma$ and consequently,
on $X^*_+$. In other words, $x\in X_+$, which we just have
seen to be wrong. So, we may  normalize $x$, such that
$\inf_{\zeta\in\Pi}\langle x,\zeta\rangle=-1$ and
\begin{equation}
\langle x,f\rangle>\lambda
\end{equation}
by (2.5).  Noting, that
$-\Pi^\circ\subset -(V-X_+)^{\circ\circ}=\cl(V+X_+)$, we get
\begin{equation}
x\in-\Pi^\circ\cap\cl C\subset\cl(V+X_+)\cap\cl C\subset
   \cl C_V.
\end{equation}
To prove the last inclusion in (2.7) note, that $\alpha x$ is an
interior point of $V+X_+$ for all $\alpha\in [0,1)$, see e.g.
\cite[Chap.\,II]{S66}. For fixed $0\le\alpha<1$ let $W$ be
a $\tau$-neighbourhood of zero such that $\alpha x+W\subset V+X_+$.
Since $\alpha x\in\cl C$, the set $ (\alpha x+W)\cap C$ is non-empty.
By (2.3) it means that $\alpha x\in \cl C_V$ for each $0\le\alpha<1$
and therefore also for $\alpha=1$.

Clearly, relations (2.6), (2.7) yield the desired contradiction to
(2.4), which completes the proof. \end{proof}

The conditions of theorem 1 are satisfied for any Banach lattice $X$ (with
the norm topology $\tau$) since we can take $\Gamma=B_{X^*}\cap X^*_+$,
where $B_{X^*}$ is the unit ball of $X^*$.
Moreover, in this case, we can consider only one neighbourhood
of zero $V=B_X$ in condition (i). The corresponding result
for the space $L^\infty$  with the norm topology is formulated
below.
\begin{cor} For any element $f\in (L^\infty)^*$ the following
conditions are equivalent:
\item[(i)]
$\sup_{x\in C_1}\langle x,f\rangle<+\infty,\ \ \
   C_1=\{x\in C: x^-\le 1\ a.s.\};$
\item[(ii)] there exists $g\in  (L^\infty)^*$
such that $g\ge f$ and $g\in C^\circ$.
\end{cor}

As a second example, the Mackey topology $\tau(L^\infty,L^1)$ is
locally convex-solid, see \cite[section 11]{AB85},
and the set
$$\Gamma=\{x\in L^\infty_+:\|x\|_{L^\infty}\le 1\}\subset L^1_+$$
is $\sigma (L^1,L^\infty)$-compact (weakly compact in $L^1$).
Thus, theorem 1 is valid for the space
$(L^\infty, \tau(L^\infty,L^1))$. To make this result
more concrete, we remind another descriptions of the topology
$\tau(L^\infty,L^1)$.

A function $\varphi: [0,\infty)\mapsto [0,\infty)$ is called
$N$-function if it is convex  and
$$ \lim_{t\to +0}\frac{\varphi(t)}{t}=0,\ \ \
   \lim_{t\to +\infty}\frac{\varphi(t)}{t}=\infty.$$
It follows that $\varphi$ is non-decreasing and continuous.
Let $\|x\|_\varphi$ denote the Luxemburg norm (see e.g.
\cite{KR61}):
$$\|x\|_\varphi=\inf\{\lambda>0:\int_\Omega\varphi(|x|/\lambda)
                 \,d\mathbf P\le 1\}.$$
It is known, that the Mackey topology $\tau(L^\infty,L^1)$ is
generated by the family of Luxemburg norms $\{\|\cdot\|_\varphi:
\varphi\in\Phi_N\}$, where $\Phi_N$ is the collection of all
$N$-functions (see \cite{N90}).

In addition, this topology is generated by sets
$$ \mu \bigcap_{k=1}^\infty U_{\varepsilon_k},\ \
   U_{\varepsilon_k}=\{x:\mathbf P(|x|\ge k)\le\varepsilon_k\},\ \
  k=1,\dots,\infty,\ \ \mu>0,$$
where $(\varepsilon_k)_{k=1}^\infty$ is any positive sequence.
Indeed, for any sequence $\varepsilon_k>0$ there exists
$N$-function $\varphi$, satisfying the conditions
$$ \varphi(t)\ge \max_{1\le i\le k}\{1/\varepsilon_i\},\ \ \
  t\ge k.$$
If $\|x\|_\varphi\le 1$ then
$$\mathbf P(|x|\ge k)=\int_{\{|x|\ge k\}}\,d\mathbf P\le
  \varepsilon_k\int_{\{|x|\ge k\}}\varphi(|x|)\,d\mathbf P
  \le\varepsilon_k.$$

Conversily, for any $N$-function $\varphi$ put
$\varepsilon_k=k^{-2}/\varphi(k+1).$
If $x\in\cap_{k=1}^\infty U_{\varepsilon_k}$, then
\begin{eqnarray*}
 \|x\|_\varphi &\le & \int_{|x|<1}\varphi(|x|)\,d\mathbf P+
 \sum_{k=1}^\infty\int_{k\le |x|<k+1}\varphi(|x|)\,d\mathbf P\\
 &\le & \varphi(1)+\sum_{k=1}^\infty\varphi(k+1)\mathbf P\{|x|\ge k\}\le
 \varphi(1)+\sum_{k=1}^\infty k^{-2}.
\end{eqnarray*}

We collect these results in the following corollary, giving
the answer to the question {\bf (Q)}.
\begin{cor}
For any element $f\in L^1$ the following conditions
are equivalent:
\item[(i)] there exists a sequence $\varepsilon_k>0$
such that
$$\sup\{\langle x,f\rangle:
   x\in\cap_{k=1}^\infty C^{\varepsilon_k}\}<\infty,\ \
   C^{\varepsilon_k}=\{x\in C: \mathbf P(x^-\ge k)\le\varepsilon_k\};$$
\item[(ii)] there exists $N$-function $\varphi$ such that
$$ \sup_{x\in C_\varphi}\langle x,f\rangle<\infty, \ \
   C_\varphi=\{x\in C:\|x^-\|_\varphi\le 1\};$$
\item[(iii)] there exists a convex solid
$\tau(L^\infty,L^1)$-neighbourhood of zero $V$ such that
$$\sup_{x\in C_V}\langle x,f\rangle<+\infty,\ \ \
   C_V=\{x\in C: x^-\in V\};$$
\item[(iv)] there exists $g\in L^1$ such that $g\ge f$ and
$g\in C^\circ$ .
\end{cor}

The equivalence between (iii) and (iv) follows from theorem 1. The
two other equivalencies are implied by the properties of the Mackey
topology $\tau(L^\infty,L^1)$, presented above.

\section{Examples}
\setcounter{equation}{0}
Recall that $(L^\infty)^*$ may be identified with the space
of all bounded finitely additive measures $\mu$ on $\mathcal F$
with the property that $\mathbf P(A)=0$ implies that $\mu(A)=0$
\cite{DunS58}. Our first example shows that in the context
of Corollary 1, in general, it is not possible to find the element
$g\in (L^\infty)^*$ already in $L^1$ even if $f\in L^\infty$.

\vspace{0.7pc}
{\bf Example 1.} Let $\Omega=[0,1]$, $\mathcal F$ consists of
all Lebesgue measurable sets and let $\mathbf P$ be the Lebesgue measure.
Consider a purely finitely additive measure $\mu:\mathcal F
\mapsto\{0,1\}$ such $\mu(I)=1$ for any open interval $I\subset (0,1)$,
containing $1/2$ (see \cite{YH52}). It follows
that $\mu\{|t-1/2|\ge\delta\}=0$ for all $\delta>0$. Put
$$C=\{x\in L^\infty: \int_\Omega x\, d(\mathbf P+\mu)\le 0\}.$$

The element $f=1\in (L^\infty)^*\cap L^\infty$ is bounded on the set
$C_1$, defined in Corollary 1:
$$ \langle x,1\rangle=\int_\Omega x\,d\mathbf P\le
   -\int_\Omega x\,d\mu\le 1,\ \ x\in C_1$$
and it is dominated by the element of $C^\circ\subset (L^\infty)^*$,
corresponding to the measure $\mathbf P+\mu$. However, $f$
is unbounded on any set $\cap_{k=1}^\infty C^{\varepsilon_k}$,
defined in Corollary 2(i).

To show this,  consider a sequence
$x_n\in L^\infty$, defined by the formulas
$$ x_n(t)=n,\ \ |t-1/2|\ge\varepsilon_n/2,\ \ \
   x_n(t)=-n,\ \ |t-1/2|<\varepsilon_n/2,$$
$n\ge 1$, $t\in [0,1]$. Without loss of generality, we may assume that
$\varepsilon_k>0$ monotonically tends to $0$.
Evidently, $x_n\in
\cap_{k=1}^\infty C^{\varepsilon_k}$:
$$ \int_\Omega x_n\,d(\mathbf P+\mu)=\int_0^1 x_n(t)\,dt
-n=-2n\varepsilon_n\le 0,$$
$$\mathbf P(x_n^-\ge k)=0,\ \ n<k;\ \
  \mathbf P(x_n^-\ge k)=\varepsilon_n
  \le\varepsilon_k,\ \ n\ge k.$$
But
$$ \langle x_n,1\rangle=\int_0^1 x_n(t)\,dt=n(1-2\varepsilon_n)
\to+\infty,\ \ n\to\infty.$$

Hence, by Corollary 2, $f=1$ cannot be dominated by any element
of $C^\circ\cap L^1$.
\vspace{0.7pc}

The next examples are in more financial spirit. Note, that
in both of them the cone $C$ is a subspace. This is not
substantial: passing to $C-L^\infty_+$, the results still
hold true.
\vspace{0.7pc}

{\bf Example 2.}  We consider a slight modification of an
example, given in \cite[Remark 5.5.2]{DS05}.
Let $\Omega=\mathbb N$, the sigma-algebra
$\mathcal F_0$ is generated by the sets $(\{2n-1,2n\})_{n=1}^\infty$,
and $\mathcal F=\mathcal F_1$ to be the power set of $\Omega$. Define
the probability measure $\mathbf P$ on $\mathcal F$ by
$\mathbf P\{2n-1\}= \mathbf P\{2n\}=2^{-n-1}$.
Let the asset prices $(S_t)_{t=0}^1$
at times $0$ and $1$ be $S_0\equiv 0$,
$$ S_1(2n-1)=1, \ \ S_1(2n)=-2^{-n},\ \ n\in\mathbb N.$$

Let the cone $C$ be generated by the elements $\gamma (S_1-S_0)$
in $L^\infty$, where $\gamma$ is $\mathcal F_0$-measurable random
variable. As usual, $\gamma$ may be interpreted as investor's
portfolio at time $t=0$. Then the set $C$ consists of possible
investor's gains at time $t=1$.
Evidently, the no-arbitrage condition (1.1) is satisfied.

We claim that for any $f\in L^1_+$ the conditions of
Corollaries 1 and 2 are equivalent and there exists an element
$g\ge f$, $g\in C^\circ\cap L^1$ if and only if
\begin{equation}
\sum_{n=0}^\infty f(2n-1)<\infty.
\end{equation}

It suffices to show that condition (3.1) implies condition
(iv) of Corollary 2 and that condition (i) of Corollary 1
implies (3.1).
Assume that (3.1) is satisfied and put
$$ g(2n-1)=\max\{f(2n-1), 2^{-n} f(2n)\},\ \
   g(2n)=2^n g(2n-1),\ \ n\in\mathbb N.$$
Then $g\in L^1(\mathbf P)$ and $g\ge f$. Computing the
conditional expectation:
$$\mathbf E_\mathbf P (g S_1|\mathcal F_0)(2n-1)=
(g(2n)S_1(2n)+g(2n-1)S_1(2n-1))/2^{n+1}=0,$$
we see that $g\in C^\circ$.

Now assume that condition (i) of Corollary 1
is satisfied. Put $\gamma(2n-1)=\gamma(2n)=2^n$. Then
$\gamma S_1\in C_1$ and
$$\langle \gamma S_1,f\rangle=\sum_{n=1}^\infty
(f(2n-1)/2-2^{-n-1}f(2n))<+\infty.$$
Since $f\in L^1(\mathbf P)$ we have
$\sum_{n=1}^\infty 2^{-n-1}f(2n)<+\infty$ and
the condition (3.1) holds true.
\vspace{0.7pc}

For the cone, considered in example 2, there is no difference
between the conditions of Corollaries 1 and 2 (in contrast to
example 1, which did not allow for a financial interpretation).
Below we consider a market with infinitely many assets,
where these conditions are different and the following is true:
\begin{equation}
(f+L_+^1)\cap C^\circ=\emptyset,\ \ \
(f+(L^\infty)^*_+)\cap C^\circ\neq\emptyset
\end{equation}
for some $f\in L^1_+$.
\vspace{0.7pc}

{\bf Example 3.} Consider the probability space $(\Omega,\mathcal F,
\mathbf P)$ as in example 1.
Let $(A_n)_{n=1}^\infty$, $A_n\subset [0,1/2]$
be a sequence of independent events with
probabilities $\mathbf P(A_n)=1/2^n$. To construct such
sequence take independent random variables $\xi_n:\Omega\mapsto\{0,1\}$
such that $\mathbf P(\xi_n=1)=1/2^{n-1}$ and put
$$ A_n=\{\xi_n^{-1}(1)\}/2=\{t\in [0,1/2]:\xi_n(2t)=1\}.$$

Furthemore, put $b_0=1/2$, $b_n=b_{n-1}+4^{-n}$, $n\ge 1$ and
consider the sequence of intervals $B_n=(b_{n-1},b_n]\subset (1/2,5/6]$.
The sets $B_n$ are mutually disjoint and disjoint from
$\cup_{n=1}^\infty A_n$. Let
$$f=\sum_{n=1}^\infty 2^n I_{B_n}+I_{[0,1/2]}+I_{[5/6,1]}.$$
Clearly, $f\in L^1_+(\mathbf P)$.

Now we introduce a countable sequence of asset price increments:
$$x_n=S_1^n-S_0^n=2^n I_{B_n}-I_{A_n}, \ \ \ n\in\mathbb N$$
at times $0$ and $1$. We assume that the processes
$(S_t^n)_{t=0}^1$ are adapted to the
filtration $(\mathcal F_0, \mathcal F_1)$, where
$\mathcal F_1=\mathcal F$ and $\mathcal F_0$ is trivial.
Portfolios $\gamma^n$ are non-random, since they are assumed
to be $\mathcal F_0$-measurable.

Let $C$ be the linear subspace of $L^\infty$ spanned (algebraically)
by $x_n$. Elements of $C$ describe the investor's gains,
obtained by trading in a finite collection of assets.
Condition  $\mathbf E_\mathbf P(x_n)=0$ imply that $C$ is disjoint from
$L^\infty_+\backslash\{0\}$.

Let $z=\sum_{n\in J}\gamma^n x_n$ be any element of $C_1$.
Here $J$ is a finite subset of $\mathbb N$ and $\gamma^n$ are some
constants. By definition of $C_1$ we have
$$z=\sum_{n\in J}\gamma^n (2^n I_{B_n}-I_{A_n}) \ge -1,\
a.s.$$
Considering this inequality on the sets $B_n$ and
$\cap_{n\in J} A_n$, we get
$$-\gamma^n 2^n\le 1, \ \ \ \sum_{n\in J}\gamma^n\le 1.$$
It follows that condition (i) of Corollary 1 is satisfied:
\begin{eqnarray*}
\langle z,f\rangle &=& \sum_{n\in J}\gamma^n(2^n\int_{B_n} f\,
d\mathbf P-\int_{A_n}f\,d\mathbf P)\\
&=&\sum_{n\in J}\gamma^n(1-2^{-n})
\le 1+\sum_{n\in J} 4^{-n}\le 4/3.
\end{eqnarray*}

To show that condition (i) of Corollary 2 fails, consider any
sequence $\varepsilon_k>0$, $k\ge 1$ and assume that
$f$ is bounded from above by a constant $\beta$ on the
set $\cap_{k=1}^\infty C_{\varepsilon_k}$. Define natural numbers
$m$, $n_1,\dots, n_m$ as follows:
$$ m>\beta+1,\ \ \ \sum_{i=1}^m\frac{1}{2^{n_i}}
   \le\min \{\varepsilon_1,\dots,\varepsilon_m\}.$$
We have $\mathbf P(x_{n_1}+\dots+x_{n_m}\le -k)=0$, $k>m$ and
$$\mathbf P(x_{n_1}+\dots+x_{n_m}\le -k)\le\mathbf P(\cup_{i=1}^m
  \{x_{n_i}\le -1\})\le
   \sum_{i=1}^m\frac{1}{2^{n_i}}\le\varepsilon_k,\ \ k\le m.$$
Thus $x_{n_1}+\dots+x_{n_m}\in\cap_{k=1}^\infty C_{\varepsilon_k}$ and
we obtain a contradiction:
\begin{eqnarray*}
   \langle x_{n_1}+\dots+x_{n_m},f\rangle &=&\sum_{i=1}^m
   \left(2^{n_i}\int_{B_{n_i} } f\,d\mathbf P-\int_{A_{n_i}}f\,d\mathbf P
   \right)\\
   &=& m-\sum_{i=1}^n 2^{-n_i}\ge m-1>\beta.
\end{eqnarray*}

Note also, that if $\nu$ is the non-negative finitely additive measure,
corresponding to an element $g\in C^\circ$, $g\ge f$, then
$$\nu(A_n)=\langle I_{A_n},g\rangle=
2^n\langle I_{B_n},g\rangle\ge 2^n\langle I_{B_n},f\rangle=1.
$$
Hence, $\nu$ is not countably additive.
\vspace{0.7pc}

 Finally, we mention that it would be interesting to clarify if
the relations (3.2) can hold true for the case of finitely many
assets.

\end{document}